\documentclass[12pt,smallextended]{svjour3}
\smartqed
\usepackage{lineno}
\modulolinenumbers[5]

\journalname{Journal of XXX}
\usepackage{amssymb,amsmath,amsfonts}
\usepackage[bookmarksnumbered, colorlinks, urlcolor=blue, linkcolor=blue, citecolor=red, plainpages=false, pdfstartview=FitW]{hyperref}
\usepackage{array}
\usepackage{bm}
\usepackage{epsfig}
\usepackage{mathptmx}
\usepackage{diagbox}
\usepackage{graphicx}
\usepackage{color}
\usepackage{caption}
\usepackage{wrapfig}
\usepackage{graphics}
\usepackage{enumerate}
\usepackage{amsxtra}
\usepackage{latexsym}
\usepackage{multirow}
\usepackage{slashbox}
\usepackage{subfigure}
\usepackage{epstopdf}
\usepackage{pdfpages}
\usepackage{algorithm}
\usepackage{algorithmic}
\usepackage{url}
\newcommand{\thmlist}{
\begin{list}{Step 1}
{\setlength{\leftmargin}{0.6 in}\setlength{\labelwidth} {0.5 in}}}

\newcommand{\alglist}{
\begin{list}{Step 1}
{\setlength{\leftmargin}{1.1 in} \setlength{\labelwidth}{1.0 in}}}

 \renewcommand{\proof} {\noindent {\bf Proof.} \quad}

\renewcommand{\subtitle}[1]{\color{blue}}

\def\red#1{\color{red}{#1}\color{black}}

\begin{document}

%\begin{frontmatter}

\title{Explicit pseudo-transient continuation and 
the trust-region updating strategy for unconstrained optimization}
\titlerunning{Explicit pseudo-transient continuation and
the trust-region updating strategy}
\author{Xin-long Luo\textsuperscript{$\ast$}  \and Hang Xiao \and Jia-hui Lv \and Sen Zhang}
\authorrunning{Luo, Xiao, Lv and Zhang}

\institute{Xin-long Luo
             \at
              Corresponding author. School of Artificial Intelligence, \\
Beijing University of Posts and Telecommunications, P. O. Box 101, \\
Xitucheng Road  No. 10, Haidian District, 100876, Beijing China\\
             \email{luoxinlong@bupt.edu.cn}            %  \\
           \and
           Hang Xiao
     \at
     School of Artificial Intelligence, \\
     Beijing University of Posts and Telecommunications, P. O. Box 101, \\
     Xitucheng Road  No. 10, Haidian District, 100876, Beijing China \\
     \email{xiaohang0210@bupt.edu.cn}
     \and
        Jia-hui Lv \at
        School of Artificial Intelligence, \\
Beijing University of Posts and Telecommunications, P. O. Box 101, \\
Xitucheng Road  No. 10, Haidian District, 100876, Beijing China\\
             \email{jhlv@bupt.edu.cn}
       \and
        Sen Zhang \at
        School of Artificial Intelligence, \\
Beijing University of Posts and Telecommunications, P. O. Box 101, \\
Xitucheng Road  No. 10, Haidian District, 100876, Beijing China\\
             \email{senzhang@bupt.edu.cn}
}

\date{Received: date / Accepted: date}
\maketitle

\begin{abstract}
This paper considers an explicit continuation method and the trust-region
updating strategy for the unconstrained optimization problem. Moreover, in 
order to improve its computational efficiency and robustness, the new method 
uses the switching preconditioning technique. In the well-conditioned phase, 
the new method uses the L-BFGS method as the preconditioning technique in 
order to improve its computational efficiency. Otherwise, the new method uses 
the inverse of the Hessian matrix as the pre-conditioner in order to improve
its robustness. Numerical results aslo show that the new method is more robust 
and faster than the traditional optimization method such as the trust-region 
method and the line search method. The computational time of the new method is 
about one percent of that of the trust-region method (the subroutine fminunc.m 
of the MATLAB2019a environment, it is set by the trust-region method) or one 
fifth of that the line search method (fminunc.m is set by the quasi-Newton method) 
for the large-scale problem. Finally, the global convergence analysis of the new 
method is also given. 
\end{abstract}



\vskip 2mm

%\textbf{AMS subject classifications.} 65H17 \and 65J15 \and 65K05 \and 65L05

%\end{frontmatter}

% \linenumbers
% main text

\section{Introduction} \label{SUBINT}

\vskip 2mm

In this article, we consider the following unconstrained optimization problem
\begin{align}
  &\min_{x \in \Re^n} \; f(x), \label{UNOPT}
\end{align}
where $f: \; \Re^{n} \to \Re$ is a differentiable function. For this problem,
there are many efficient methods to solve it such as the line search method and
the trust-region method \cite{CGT2000,NW1999,SY2006,Yuan2015}. The continuation
method \cite{AG2003,Botsaris1978,Goh2010,HM1996,KK1998,LQQ2004,Tanabe1980}
is another method other than the traditional optimization methods for the problem
\eqref{UNOPT}. The advantage of the continuation method over the line search method
or the trust-region method is capable of finding many local optimal points
of the non-convex optimization problem by tracking its trajectory, and it is
even possible to find the global optimal solution
\cite{BB1989,Schropp1995,Yamashita1980}. However, the computational time of
the continuation method  may be higher than that of the traditional
optimization method.

\vskip 2mm

Recently, Luo, Xiao and Lv \cite{LXL2020} consider the continuation Newton method
with the adaptive time-stepping size based on the trust-region updating strategy
for the nonlinear system of equations. According to their numerical experiments,
their method is robust and efficient to solve the nonlinear system of equations.
In order to improve the computational efficiency and the robustness of their
continuation method \cite{LXL2020} further for the large-scale optimization
problem, we consider the switching preconditioning technique for the ill-conditioned
problem. That is to say, in the well-conditioned phase, the new method uses 
the L-BFGS method as the  preconditioning technique in order to improve its 
computational efficiency. Otherwise, the new  method uses the inverse of 
the Hessian matrix as the pre-conditioner in order to improve its robustness.

\vskip 2mm

The rest of the paper is organized as follows. In section 2, we give a new
continuation method with the trust-region updating strategy and the
switching preconditioning technique for the unconstrained optimization problem
\eqref{UNOPT}. In section 3, we analyze the global convergence of this new
method. In section 4, we report some promising numerical results of the new
method, in comparison to the traditional trust-region method (fminunc.m is 
set by the trust-region method) and the line search method (fminunc.m is set 
by the quasi-Newton method) for some large-scale problems. Finally, we give 
some discussions and conclusions in section 5.

\vskip 2mm

\section{Explicit pseudo-transient continuation}

\vskip 2mm

In this section, we construct an explicit continuation method with the adaptive 
time-stepping scheme based on the trust-region updating strategy \cite{Yuan2015}
for the unconstrained optimization problem \eqref{UNOPT}. Firstly, we construct 
a generalized gradient flow for the stationary point of the unconstrained 
optimization \eqref{UNOPT}. Then, we construct an explicit continuation method 
with an adaptive time-stepping scheme for this special ordinary differential 
equations (ODEs). In order to improve the robustness of the new method, we also 
consider a switching preconditioning technique between the L-BFGS method and the 
inverse of the Hessian matrix.

\vskip 2mm

\subsection{The generalized gradient flow}

\vskip 2mm

For the unconstrained optimization problem \eqref{UNOPT}, we consider the damped
Newton method \cite{SY2006} as follows:
\begin{align}
      x_{k+1} = x_{k}  - \alpha_{k} B(x_{k})^{-1} g(x),
      \label{DNEWTON}
\end{align}
where $B(x) = \nabla^{2} f(x)$ and $g(x) = \nabla f(x)$. If we regard
$x_{k} = x(t_{k})$, $x_{k+1} = x(t_{k}+ \alpha_{k})$ and let
$\alpha_{k} \to 0$, we obtain the continuous Newton flow
\cite{Branin1972,Davidenko1953,Tanabe1979} as follows:
\begin{align}
    \frac{dx}{dt} = - B(x)^{-1} g(x), \;  x(0) = x_{0}. \label{NEWTONFLOW}
\end{align}
Actually, if we apply an iteration with the explicit Euler method
\cite{SGT2003} for the continuous Newton flow \eqref{NEWTONFLOW}, we
obtain the damped Newton method \eqref{DNEWTON}. Since $B(x)$ may be singular, we
reformulate the continuous Newton flow \eqref{NEWTONFLOW} as the following
general formula:
\begin{align}
     B(x)\frac{dx}{dt} + g(x) = 0, \; x(0) = x_{0}. \label{DAEFLOW}
\end{align}
For the continuous Newton flow \eqref{DAEFLOW}, we have the following property
\ref{PRODAEFLOW}.

\vskip 2mm

\begin{property} (Branin \cite{Branin1972} and Tanabe \cite{Tanabe1979})
\label{PRODAEFLOW} Assume that  $x(t)$ is the solution of the continuous Newton flow
\eqref{DAEFLOW}, then $r(x(t)) = \|g(x(t))\|^{2}$ converges to zero when
$t \to \infty$. That is to say, for every limit point $x^{\ast}$ of $x(t)$, it is
also a stationary point of the Newton flow \eqref{DAEFLOW} and  every element
$g_{i}(x(t))$ of $g(x(t))$ has the same convergence rate $e^{-t}$. Furthermore, 
$x(t)$ can not converge to the stationary point $x^{\ast}$ of the continuous Newton 
flow \eqref{DAEFLOW} on the finite interval when the initial point $x_{0}$ is not 
its stationary point.
\end{property}
\proof Assume that $x(t)$ is the solution of the continuous Newton flow
\eqref{DAEFLOW}, then we have
\begin{align}
    \frac{d}{dt} \left(e^{t}g(x) \right) = e^{t} B(x) \frac{dx(t)}{dt}
    + e^{t} g(x) = 0. \nonumber
\end{align}
Consequently, we obtain
\begin{align}
     g(x(t))  = g(x_{0})e^{-t}.    \label{FUNPAR}
\end{align}
From equation \eqref{FUNPAR}, it is not difficult to know that every element
$g_{i}(x(t))$ of $g(x(t))$ converges to zero with the linear convergence
rate $e^{-t}$ when $t \to \infty$. Thus, if the solution $x(t)$ of the continuous
Newton flow \eqref{DAEFLOW} belongs to a compact set, it has a limit point
$x^{\ast}$ when $t \to \infty$, and this limit point $x^{\ast}$ is
also its stationary point.

\vskip 2mm

If we assume that the solution $x(t)$ of the continuous Newton flow
\eqref{DAEFLOW} converges to its stationary point $x^{\ast}$ on the finite interval
$(0, \, T]$, from equation \eqref{FUNPAR}, we have
\begin{align}
      g(x^{\ast}) = g(x_{0})e^{-T}.    \label{FLIMT}
\end{align}
Since $x^{\ast}$ is a stationary point of the continuous Newton flow \eqref{DAEFLOW},
we have $g(x^{\ast}) = 0$. By substituting it into equation \eqref{FLIMT}, we
obtain
\begin{align}
    g(x_{0}) = 0. \nonumber
\end{align}
Thus, it contradicts the assumption that $x_{0}$ is not a stationary point of the
continuous Newton flow \eqref{DAEFLOW}. Therefore, the solution $x(t)$ of
the continuous Newton flow \eqref{DAEFLOW} can not converge to its stable point
$x^{\ast}$  on the finite interval. \qed

\vskip 2mm

We can also regard the continuous Newton flow as the generalized gradient
flow (p. 361, \cite{HM1996}):
\begin{align}
    \frac{dx}{dt} = - H(x) g(x) ,
    \; x(0) = x_{0}, \label{GGF}
\end{align}
where $H(x)$ equals the inverse $B(x)^{-1}$ of the Hessian matrix  or its
quasi-Newton approximation. $H(x)$ can be regarded as a pre-conditioner of
$g(x)$ to mitigate the stiffness of the ODEs \eqref{GGF}. Consequently,
we can adopt the explicit method to compute the trajectory of the ODEs
\eqref{GGF} efficiently \cite{LXL2020}.

\vskip 2mm

\begin{remark}
If we assume that $x(t)$ is the solution of the ODEs \eqref{GGF} and $H(x)$ is
a symmetric positive definite matrix, we obtain
\begin{align}
     \frac{df(x)}{dt} = g(x)^{T} \frac{dx}{dt}
    = - g(x)^{T} H(x) g(x)  \le 0.
     \nonumber
\end{align}
That is to say, $f(x)$ is monotonically decreasing along the solution curve $x(t)$
of the dynamical system \eqref{GGF}. Furthermore, the solution $x(t)$ converges
to $x^{\ast}$ when $f(x)$ is lower bounded and $t$ tends to infinity
\cite{HM1996,LQQ2004,Schropp1995,Tanabe1980}, where $x^{\ast}$ is the stationary 
point of the generalized gradient flow \eqref{GGF}.  Thus, we can follow
the trajectory $x(t)$ of the ODEs \eqref{GGF} to obtain its stationary point
$x^{\ast}$, which is also one stationary point of the original optimization problem
\eqref{UNOPT}.
\end{remark}

\vskip 2mm

\subsection{The explicit continuation method} \label{SUBSICM}

\vskip 2mm

The solution curve $x(t)$ of the OEDs \eqref{GGF} can not be efficiently 
followed on an infinite interval by the traditional ODE method
\cite{AP1998,BJ1998,JT1995,SGT2003}, so we need to construct the particular
method for this problem \eqref{GGF}. We apply the first-order implicit Euler
method \cite{AP1998,SGT2003} to the ODEs \eqref{GGF}, then we obtain
\begin{align}
      x_{k+1} =  x_{k} - \Delta t_{k} H(x_{k+1})g(x_{k+1}),
      \label{IEGPGF}
\end{align}
where $\Delta t_k$ is the time-stepping size.

\vskip 2mm

Since the system of equations \eqref{IEGPGF} is a nonlinear system which can not
be directly solved, we seek for its explicit approximation formula.
We denote $s_{k} = x_{k+1} - x_{k}$. By using the first-order Taylor expansion,
we have the linear approximation $g(x_{k}) + B(x_{k})s_{k}$ of
$g(x_{k+1})$. By substituting it into equation \eqref{IEGPGF} and using the
zero-order approximation $H(x_{k})$ of $H(x_{k+1})$, we have
\begin{align}
    s_{k} \approx - \Delta t_{k} H(x_{k})(g(x_{k}) + B(x_{k})s_{k})
    = - \Delta t_{k} H(x_{k}) g(x_{k}) - \Delta t_{k} H(x_{k})B(x_{k})s_{k}.
    \label{AIEGPGF}
\end{align}
Let $H(x_{k}) = \left(\nabla^{2} f(x_{k})\right)^{-1}$. Then, we have $H(x_{k})B(x_{k}) = I$.
By substituting it into equation \eqref{AIEGPGF}, we obtain the explicit
continuation method as follows:
\begin{align}
      & s_{k}^{N} = - H_{k}g_{k}, \;  
      s_{k} = \frac{\Delta t_{k}}{1 + \Delta t_{k}}s_{k}^{N}, \label{EXPTC} \\
      & x_{k+1} = x_{k} + s_{k}, \label{XK1}
\end{align}
where $g_{k} = \nabla f(x_{k})$ and $H_{k} = (\nabla^{2} f(x_{k}))^{-1}$ or
its quasi-Newton approximation.

\vskip 2mm

\begin{remark} 
The explicit continuation method \eqref{EXPTC}-\eqref{XK1} equals the damped
Newton method if we let $\alpha_{k} = \Delta t_k/(1+\Delta t_k)$ in equation
\eqref{EXPTC}. However, from the view of the ODE method, they are different.
The damped Newton method \eqref{DNEWTON} is derived from the explicit Euler
method applied to the generalized gradient flow \eqref{GGF}. Its time-stepping
size $\alpha_k$ is restricted by the numerical stability \cite{SGT2003}. That
is to say, for the linear test equation $dx/dt = - \lambda x$, its time-stepping
size $\alpha_{k}$ is restricted by the stable region $|1-\lambda \alpha_{k}| \le 1$.
Therefore, the large time-stepping size $\alpha_{k}$ can not be adopted in
the steady-state phase.

\vskip 2mm

The explicit continuation method \eqref{EXPTC}-\eqref{XK1} is derived from
the implicit Euler scheme and the linear approximation $g(x_{k}) + B(x_{k})s_{k}$ 
of $g(x_{k+1})$ applied to the generalized gradient flow
\eqref{GGF}, and its time-stepping size $\Delta t_k$ is not restricted by
the numerical stability for the linear test equation. Therefore, the large
time-stepping size $\Delta t_{k}$ can be adopted in the steady-state phase,
and the explicit continuation method \eqref{EXPTC}-\eqref{XK1} mimics the Newton
method. Consequently, it has the fast convergence rate near the stationary
point $x^{\ast}$ of the generalized gradient flow \eqref{GGF}. The most of all,
the new time-stepping size $\alpha_{k} = \Delta t_{k}/(\Delta t_{k} + 1)$ is
favourable to adopt the trust-region updating strategy for adaptively adjusting
the time-stepping size $\Delta t_{k}$ such that the explicit continuation method
\eqref{EXPTC}-\eqref{XK1} accurately tracks the generalized gradient flow
\eqref{GGF} in the transient-state phase and achieves the fast convergence rate
in the steady-state phase.
\end{remark}

\vskip 2mm

\begin{remark} \label{REMPLS}
Luo, Xiao, Lv \cite{LXL2020} and Luo, Yao \cite{LY2021} have considered the 
explicit continuation method for nonlinear equations and the linear programming 
problem when $H_{k} = (g'(x_{k}))^{-1}$, respectively. Here, we consider the 
quasi-Newton approximation $H_{k}$ of $(g'(x_{k}))^{-1}$ in equation \eqref{XK1}.
\end{remark}

\vskip 2mm

\vskip 2mm

\subsection{The trust-region updating strategy}

\vskip 2mm

Another issue is how to adaptively adjust the time-stepping size $\Delta t_k$
at every iteration. We borrow the adjustment method of the trust-region radius
from the trust-region method due to its robustness and its fast convergence
rate \cite{CGT2000,Yuan2015}. When we use the trust-region updating strategy 
to adaptively adjust time-stepping size $\Delta t_{k}$ 
\cite{Higham1999,LLT2007,LKLT2009,LLS2021}, we also need to construct a local 
approximation model of $f(x)$ around $x_{k}$. Here, we adopt the following 
quadratic function as its approximation model:
\begin{align}
     q_k(x_{k} + s) = f(x_{k}) + s^{T}g_{k} + \frac{1}{2}s^{T}B_{k}s,
     \label{QOAM}
\end{align}
where $B_{k} = \nabla^{2} f(x_{k})$ or its quasi-Newton approximation.
In practical computation, we do not store the matrix $B_{k}$. Thus, we use the
explicit continuation method \eqref{EXPTC}-\eqref{XK1} and regard $H_{k} = B_{k}^{-1}$
to simplify the quadratic model $q_k(x_{k}+s_{k}) - q(x_{k})$ as follows:
\begin{align}
    m_{k}(s_{k}) = g_{k}^{T}s_{k}
    - \frac{0.5\Delta t_{k}}{1+\Delta t_{k}}g_{k}^{T}s_{k}
   = \frac{1+0.5\Delta t_{k}}{1+\Delta t_{k}} g_{k}^{T}s_{k}
   \approx q_{k}(x_{k}+s_{k}) - q_{k}(x_{k}).  \label{LOAM}
\end{align}
where $g_{k} = \nabla f(x_k)$. We enlarge or reduce the time-stepping size
$\Delta t_k$ at every iteration according to the following ratio:
\begin{align}
    \rho_k = \frac{f(x_k)-f(x_{k+1})}{m_k(0)-m_k(s_{k})}.
    \label{MRHOK}
\end{align}
A particular adjustment strategy is given as follows:
\begin{align}
     \Delta t_{k+1} = \begin{cases}
          \gamma_1 \Delta t_k, &{if \hskip 1mm 0 \leq \left|1- \rho_k \right| \le \eta_1,}\\
          \Delta t_k, &{if \hskip 1mm \eta_1 < \left|1 - \rho_k \right| < \eta_2,}\\
          \gamma_2 \Delta t_k, &{if \hskip 1mm \left|1-\rho_k \right| \geq \eta_2,}
                   \end{cases} \label{ADTK1}
\end{align}
where the constants are selected as $\eta_1 = 0.25, \; \gamma_1 = 2, \; \eta_2 = 0.75, \;
\gamma_2 = 0.5$  according to numerical experiments. When $\rho_{k} \ge \eta_{a}$,
we accept the trial step $s_{k}$ and let $x_{k+1} = x_{k} + s_{k}$, where
$\eta_{a}$ is a small positive number such as $\eta_{a} = 1.0\times 10^{-6}$.
Otherwise, we discard it and let $x_{k+1} = x_{k}$.

\vskip 2mm

\begin{remark}
This new time-stepping size selection based on the trust-region strategy
has some advantages compared to the traditional line search strategy.
If we use the line search strategy and the damped Newton method \eqref{DNEWTON}
to track the trajectory $x(t)$ of the continuous Newton flow \eqref{DAEFLOW},
in order to achieve the fast convergence rate in the steady-state phase,
the time-stepping size $\alpha_{k}$ of the damped Newton method is tried from 1
and reduced by the half with many times at every iteration. Since the linear model
$f(x_{k}) + g_{k}^{T}s_{k}$ may not approximate $f(x_{k}+s_{k})$ well in the
transient-state phase, the time-stepping size $\alpha_{k}$ will be small.
Consequently, the line search strategy consumes the unnecessary trial steps in
the transient-state phase. However, the selection of the time-stepping size
$\Delta t_{k}$ based on the trust-region strategy \eqref{MRHOK}-\eqref{ADTK1}
can overcome this shortcoming.
\end{remark}

\subsection{The switching  preconditioning technique}

\vskip 2mm

For the large-scale problem, the numerical evaluation of the Hessian matrix
$\nabla^{2}f(x_{k})$ consumes much time. In order to overcome this shortcoming,
we use the limited-memory BFGS quasi-Newton formula
(see \cite{Broyden1970,Fletcher1970,Goldfarb1970,Mascarenhas2004,Shanno1970}
or pp. 222-230, \cite{NW1999}) to approximate $(\nabla^{2}f(x_{k}))^{-1}$.

\vskip 2mm 

Recently, Ullah, Sabi and Shah \cite{USS2020} give an efficient L-BFGS updating
formula for the system of monotone nonlinear equations. In order to avoid the
ill-conditioning of $B_{k}$, they considered the revised
BFGS updating formula \cite{USS2020} as follows:
\begin{align}
       B_{k+1} = \lambda_{k}\left(I - \frac{s_{k}s_{k}^{T}}{s_{k}^{T}s_{k}}\right)
       + \sigma_{k} \frac{y_{k}y_{k}^{T}}{y_{k}^{T}s_{k}}, \label{LBFGSP}
\end{align}
where $y_{k} = g_{k+1} - g_{k}, \; s_{k} = x_{k+1} - x_{k}$ and $\lambda_{k}, \,
\sigma_{k}$ are two undetermined parameters. Then, they solved the minimizer of
the measurement function $\varphi$ \cite{BNY1987} on the variables $\lambda_{k}$ 
and $\sigma_{k}$, where $\varphi(B_{k+1})$ is defined by
\begin{align}
      & \varphi(\lambda_{k}, \, \sigma_{k}) 
      = \Psi(B_{k+1}) = \text{trace}(B_{k+1}) - \ln(\det(B_{k+1})) \nonumber \\
      & \quad = (n-1) \lambda_{k} + \sigma_{k}\frac{\|y_{k}\|^{2}}{y_{k}^{T}s_{k}}
      - \ln\left(\lambda_{k}^{n-1}\sigma_{k} \frac{y_{k}^{T}s_{k}}{\|s_{k}\|^{2}}\right)
      \nonumber \\
      & = (n-1)(\lambda_{k}-\ln(\lambda_{k})) + \sigma_{k}\frac{\|y_{k}\|^{2}}{y_{k}^{T}s_{k}}
      - \ln(\sigma_{k}) - \ln(y_{k}^{T}s_{k}) + \ln(\|s_{k}\|^{2}). \label{VARFUN}
\end{align}
Consequently, by solving $\min \varphi(\lambda_{k}, \, \sigma_{k})$,
they obtained the optimal parameters $\lambda_{k} = 1$ and
$\sigma_{k} = y_{k}^{T}s_{k}/\|y_{k}\|^{2}$. By substituting them into equation
\eqref{LBFGSP}, they obtained the revised L-BFGS updating formula:
\begin{align}
       B_{k+1} = I - \frac{s_{k}s_{k}^{T}}{s_{k}^{T}s_{k}}
       + \frac{y_{k}y_{k}^{T}}{\|y_{k}\|^{2}}. \label{LBFGSR}
\end{align}

\vskip 2mm

By using the Sherman-Morrison-Woodburg formula (P. 17, \cite{SY2006}), from
equation \eqref{LBFGSR}, we obtain the inverse of $B_{k+1}$:
\begin{align}
      H_{k+1} = B_{k+1}^{-1} = I  - \frac{y_{k}s_{k}^{T}
       + s_{k}y_{k}^{T}}{y_{k}^{T}s_{k}}
       + 2\frac{y_{k}^{T}y_{k}}{(y_{k}^{T}s_{k})^{2}} s_{k}s_{k}^{T},
      \label{ILBFGS}
\end{align}
The initial matrix $H_{0}$ can be simply selected by the inverse 
$(\nabla^{2} f(x_{0}))^{-1}$ of the Hessian matrix. From equation 
\eqref{ILBFGS}, it is not difficult to verify
\begin{align}
       H_{k+1}y_{k} = \frac{y_{k}^{T}y_{k}}{y_{k}^{T}s_{k}} s_{k}. \nonumber
\end{align}
That is to say, $H_{k+1}$ satisfies the scaling quasi-Newton property.

\vskip 2mm

The L-BFGS updating formula \eqref{ILBFGS} has some nice properties such as the
symmetric positive definite property and the positive lower bound of its eigenvalues.

\vskip 2mm

\begin{lemma} \label{LEMHLB}
Matrix $H_{k+1}$ defined by equation \eqref{ILBFGS} is symmetric positive definite
and its eigenvalues are greater than $1/2$.
\end{lemma}

\vskip 2mm

\proof (i) For any nonzero vector $z \in \Re^{n}$, from equation \eqref{ILBFGS},
we have
\begin{align}
       & z^{T}H_{k+1}z = \|z\|^{2} - 2{(z^{T}y_{k})(z^{T}s_{k})}/{y_{k}^{T}s_{k}}
       + 2(z^{T}s_{k})^{2}{\|y_{k}\|^{2}}/{(y_{k}^{T}s_{k})^{2}} \nonumber \\
       & \quad =  \left(\|z\| -
         \left|{z^{T}s_{k}}/{y_{k}^{T}s_{k}}\right|\|y_{k}\|\right)^{2}
          + 2 \|z\| \left|{z^{T}s_{k}}/{y_{k}^{T}s_{k}}\right|\|y_{k}\|
          \nonumber \\
       & \quad \quad  - 2 {(z^{T}y_{k})(z^{T}s_{k})}/{y_{k}^{T}s_{k}}
          +\|y_{k}\|^{2} (z^{T}s_{k}/y_{k}^{T}s_{k})^{2}
          \ge 0.  \label{ZTHK1Z}
\end{align}
In the last inequality of equation \eqref{ZTHK1Z}, we use the Cauchy-Schwartz
inequality $\|z^{T}y\| \le \|z\|\|y_{k}\|$ and its equality holds if only if
$z = t y_{k}$. When $z = t y_{k}$, from equation \eqref{ZTHK1Z}, we have
$z^{T}H_{k+1}z = t^{2}\|y_{k}\|^{2} = \|z\|^{2}> 0$. When $z^{T}s_{k} = 0$,
from equation \eqref{ZTHK1Z}, we also have $z^{T}H_{k+1}z = \|z\|^{2} > 0$.
Therefore, we conclude that $H_{k+1}$ is a symmetric positive definite matrix.

\vskip 2mm

(ii) It is not difficult to know that it exists at least $n-2$ linearly independent
vectors $z_{1}, \, z_{2}, \, \ldots, \, z_{n-2}$ such that $z_{i}^{T}s_{k} = 0, \,
z_{i}^{T}y_{k} = 0 \, ( i = 1, \, 2, \, \ldots, \, (n-2))$ hold. That is to say, 
matrix $H_{k+1}$ defined by equation \eqref{ILBFGS} has at least $(n-2)$ linearly 
independent eigenvectors whose corresponding eigenvalues are 1. We denote the other 
two eigenvalues of $H_{k+1}$ as $\mu_{i}^{k+1} \, (i = 1, \, 2)$ and their 
corresponding eigenvalues as $p_{1}$ and $p_{2}$, respectively. Then, from equation 
\eqref{ILBFGS}, we know that the eigenvectors $p_{i} \, (i = 1, \, 2)$ can be 
represented as $p_{i} = y_{k} + \beta_{i} s_{k}$ when 
$\mu_{i}^{k+1} \neq 1 \, (i = 1, \, 2)$. From equation 
\eqref{ILBFGS} and $H_{k+1}p_{i} = \mu_{i}^{k+1} p_{i} \, (i = 1, \, 2)$, we have
\begin{align}
     - \left(\mu_{i}^{k+1} + \beta_{i}\frac{s_{k}^{T}s_{k}}{s_{k}^{T}y_{k}}\right)y_{k}
     + \left(\frac{y_{k}^{T}y_{k}}{y_{k}^{T}s_{k}}
     + 2 \beta_{i} \frac{(y_{k}^{T}y_{k})(s_{k}^{T}s_{k})}{(y_{k}^{T}s_{k})^{2}}
     - \mu_{i}^{k+1} \beta_{i}\right)s_{k} = 0. \label{EIGVAS}
\end{align}

\vskip 2mm

When $y_{k} = ts_{k}$, from equation \eqref{ILBFGS}, we have $H_{k+1} = I$. In this
case, we conclude that the eigenvalues of $H_{k+1}$ are greater than $1/2$. When
vectors $y_{k}$ and $s_{k}$ are linearly independent, from equation \eqref{EIGVAS},
we have
\begin{align}
      & \mu_{i}^{k+1} + \beta_{i} {s_{k}^{T}s_{k}}/{s_{k}^{T}y_{k}} = 0,
      \nonumber \\
      &
      {y_{k}^{T}y_{k}}/{y_{k}^{T}s_{k}}
      + 2 \beta_{i}{(y_{k}^{T}y_{k})(s_{k}^{T}s_{k})}/{(y_{k}^{T}s_{k})^{2}}
     - \mu_{i}^{k+1}\beta_{i} = 0, \; i = 1, \, 2. \nonumber
\end{align}
That is to say, $\mu_{i}^{k+1} \, (i = 1:2)$ are the two solutions of the following
equation:
\begin{align}
      \mu^{2} - 2\mu (y_{k}^{T}y_{k})(s_{k}^{T}s_{k})/(s_{k}^{T}y_{k})^{2}
      + (y_{k}^{T}y_{k})(s_{k}^{T}s_{k})/(s_{k}^{T}y_{k})^{2} = 0. \label{QUDEDQ}
\end{align}
Consequently, from equation \eqref{QUDEDQ}, we obtain
\begin{align}
      \mu_{1}^{k+1} + \mu_{2}^{k+1} = 2 (y_{k}^{T}y_{k})(s_{k}^{T}s_{k})/(s_{k}^{T}y_{k})^{2},
      \;
      \mu_{1}^{k+1}\mu_{2}^{k+1} = (y_{k}^{T}y_{k})(s_{k}^{T}s_{k})/(s_{k}^{T}y_{k})^{2}.
      \label{ROOTQEQ}
\end{align}
From equation \eqref{ROOTQEQ}, it is not difficult to obtain
\begin{align}
       {1}/{\mu_{1}^{k+1}} + {1}/{\mu_{2}^{k+1}} = 2, \;  \mu_{i}^{k+1} > 0,
       \; i = 1, \, 2.    \label{REROOT}
\end{align}
Therefore, from equation \eqref{REROOT}, we conclude that
$\mu_{i}^{k+1} > \frac{1}{2}  \, (i = 1, \, 2)$. Consequently, the eigenvalues of $H_{k+1}$
are greater than 1/2. \qed

\vskip 2mm

According to our numerical experiments, the quasi-Newton updating method
\eqref{ILBFGS} works well for the most unconstrained optimization problems.
However, for the very ill-conditioned problem, the quasi-Newton method
\eqref{ILBFGS} will fail to accurately obtain its minimizer $x^{\ast}$.
In order to improve the robustness of the method, we use the inverse 
$(\nabla^{2} f(x_{k}))^{-1}$ of the Hessian matrix as the pre-conditioner of 
the gradient $g(x_{k})$ in the ill-conditioned phase. Furthermore, we identify
the ill-conditioned phase by the gap between $f(x_{k}+s_{k})$ and its approximation
$f(x_{k}) + g(x_{k})^{T}s_{k} + \frac{1}{2} s_{k}^{T}B_{k}s_{k}$ when $s_{k}$
is small. That is to say, we regard the problem as an ill-conditioned phase
when the number of $k$ that satisfies $|1-\rho_{k}| \ge \eta_{2}$
is greater than the threshold. Therefore, we give the following switching
preconditioning technique as follows:
\begin{align}
    H_{k+1} = \begin{cases}
                  I  - \frac{y_{k}s_{k}^{T}+ s_{k}y_{k}^{T}}{y_{k}^{T}s_{k}}
                  + 2\frac{y_{k}^{T}y_{k}}{(y_{k}^{T}s_{k})^{2}}s_{k}s_{k}^{T},
                  \; \text{if} \; K_{bad}\le 5 \; \text{and}
                  \; |s_{k}^{T}y_{k}| > \theta \|s_{k}\|^{2}, \\
                  \left(\nabla^{2} f(x_{k+1})\right)^{-1}, \; \text{otherwise},
              \end{cases} \label{SHBFGS}
\end{align}
where $K_{bad}$ represents the number of $k$ such that $|1-\rho_{k}| \ge \eta_{2}$
holds, the ratio $\rho_{k}$ is defined by equation \eqref{MRHOK}, and $\theta$ 
is a small positive constant such as $\theta = 10^{-6}$.

\vskip 2mm

For a real-world problem, the analytical Hessian $\nabla^{2} f(x)$ can not be
offered. Thus, in practice, we replace the Hessian matrix $\nabla^{2} f(x_{k+1})$ 
with its difference approximation as follows:
\begin{align}
    \nabla^{2} f(x_{k+1}) \approx
    \left[\frac{g(x_{k+1} + \epsilon e_{1}) - g(x_{k+1})}{\epsilon}, \,
    \ldots, \, \frac{g(x_{k+1} + \epsilon e_{n}) - g(x_{k+1})}{\epsilon}\right],
    \label{NUMHESS}
\end{align}
where $e_{i}$ represents the unit vector whose elements equal zeros except
for the $i$-th element equals 1, and the parameter $\epsilon$ can be selected as
$10^{-6}$ according to our numerical experiments.

\vskip 2mm

According to the above discussions, we give the detailed implementation of
the explicit continuation method with the trust-region updating strategy 
for the unconstrained optimization problem \eqref{UNOPT} in Algorithm
\ref{ALGEPTC}.

\begin{algorithm}
	\renewcommand{\algorithmicrequire}{\textbf{Input:}}
	\renewcommand{\algorithmicensure}{\textbf{Output:}}
    \newcommand{\algorithmicbreak}{\textbf{break}}
    \newcommand{\BREAK}{\STATE \algorithmicbreak}
	\caption{Explicit pseudo-transient continuation and the trust-region updating strategy 
    for unconstrained optimization (the Eptctr method)}
    \label{ALGEPTC}	
	\begin{algorithmic}[1]
		\REQUIRE ~~\\
        the objective function $f(x)$, the initial point $x_0$ \, (optional),
        the terminated parameter $\epsilon$ \, (optional).
		\ENSURE ~~\\
        the optimal approximation solution $x^{\ast}$.

        \vskip 2mm
        		
        \STATE Set $x_0 = 2*\text{ones}(n, \, 1)$ and $\epsilon = 10^{-6}$ as the default values.
        \STATE Initialize the parameters: $\eta_{a} = 10^{-6}, \; \eta_1 = 0.25,
        \; \gamma_1 =2, \; \eta_2 = 0.75, \; \gamma_2 = 0.5, \; \theta = 10^{-6}$.
        \STATE Set $k = 0$. Evaluate $f_0 = f(x_0), \; g_0 = \nabla f(x_0)$ and
        $B_{0} = \nabla^{2} f(x_{0})$.
        \STATE Set $y_{-1} = 0, \; s_{-1} = 0, \; K_{bad} = 0$, flag\_success\_trialstep = 1, 
        and $\Delta t_0 = 10^{-2}$.
        \STATE Solve $s_{0}^{N}B_{0} = - g_{0}$ to obtain $s_{0}^{N}$. 
        \WHILE{$\left(\|g_k\|> \epsilon\right)$}
          \IF{(flag\_success\_trialstep = == 1)}
          \IF{$\left(\left(|s_{k-1}^{T}y_{k-1}| > \theta s_{k-1}^{T}s_{k-1}\right)
          \&\& \left(K_{bad} < 5\right)\right)$}
            \STATE $s_{k}^{N} = - \left(g_{k} - \frac{y_{k-1}(s_{k-1}^{T}g_{k})
                     + s_{k-1}(y_{k-1}^{T}g_{k})}{y_{k-1}^{T}s_{k-1}}
                     + 2 \frac{\|y_{k-1}\|^{2}(s_{k-1}^{T}g_{k})}
                      {(y_{k-1}^{T}s_{k-1})^{2}}s_{k-1}\right)$.
          \ELSE
            \STATE Solve $B_{k} s_{k}^{N} = - g_{k}$ to obtain $s_{k}^{N}$.
          \ENDIF
          \ENDIF
          \STATE Compute $s_{k} = \frac{\Delta t_{k}}{1 + \Delta t_{k}}s_{k}^{N}$ and 
          $x_{k+1} = x_{k} + s_{k}$.
          \STATE Evaluate $f_{k+1} = f(x_{k+1})$ and compute the ratio $\rho_{k}$
          from equations \eqref{LOAM}-\eqref{MRHOK}.
          \IF{$\rho_k\le \eta_{a}$}
            \STATE Set $x_{k+1} = x_{k}, \; f_{k+1} = f_{k},
            \;  g_{k+1} = g_{k}, \; y_{k} = y_{k-1}$, flag\_success\_trialstep = 0.
          \ELSE
            \STATE Compute $g_{k+1} = \nabla f(x_{k+1}), \; y_{k} = g_{k+1} - g_{k}, \; 
            s_{k} = x_{k+1} - x_{k}$. Set flag\_success\_trialstep = 1.
          \ENDIF
          \IF{$|1 - \rho_{k}| \ge \eta_{2}$}
              \STATE $K_{bad} = K_{bad} + 1$; $\Delta t_{k+1} = \gamma_{2} \Delta t_{k}$.
          \ELSIF{$|1 - \rho_{k}| \ge \eta_{1}$}
              \STATE $\Delta t_{k+1} = \Delta t_{k}$.
          \ELSE
              \STATE $\Delta t_{k+1} = \gamma_{1} \Delta t_{k}$.
          \ENDIF
          \IF{$\left(\left(|s_{k-1}^{T}y_{k-1}| \le \theta s_{k-1}^{T}s_{k-1}\right)
          \text{or} \left(K_{bad} \ge 5\right)\right)$}
               \STATE Evaluate $B_{k+1} = \nabla^{2} f(x_{k+1})$ from the difference approximation
               \eqref{NUMHESS}.
          \ENDIF
          \STATE Set $k \leftarrow k+1$.
        \ENDWHILE
	\end{algorithmic}
\end{algorithm}

\section{Algorithm analysis}

In this section, we analyze the global convergence of the explicit continuation
method \eqref{EXPTC}-\eqref{XK1} with the trust-region updating strategy and the
switching preconditioning technique \eqref{SHBFGS} for the unconstrained optimization
problem (i.e. Algorithm \ref{ALGEPTC}). Firstly, we give a lower-bounded
estimation of $m_{k}(0) - m_{k}(s_{k})$ $(k = 1, \, 2, \, \ldots)$. This
result is similar to that of the trust-region method for the unconstrained
optimization problem \cite{Powell1975}. We denote the level set $S_{f}$ as
\begin{align}
     S_{f} = \{x: \; f(x) \le f(x_0)\}. \label{LSCONFBD}
\end{align}

\vskip 2mm

\begin{lemma} \label{LBSOAM}
Assume that there exist two positive constants $m$ and $M$ such that 
\begin{align}
     M \|z\|^{2} \ge z^{T}\nabla^{2} f(x)z \ge m \|z\|^{2} > 0, \;
     \forall x \in S_{f} \label{ASCONV}
\end{align}
holds for all $z \in \Re^{n}$. Furthermore, we assume that the quadratic model 
$q_{k}(x)$ is defined by equation \eqref{LOAM} and $s_{k}$ is computed by the 
explicit continuation method \eqref{EXPTC}-\eqref{XK1} and the switching 
preconditioning formula \eqref{SHBFGS}.
Then, we have
\begin{align}
    m_{k}(0) - m_{k}(s_{k}) \ge \frac{c_{m}\Delta t_{k}}{2(1+\Delta t_{k})}
    \left\|g_{k} \right\|^{2},     \label{PLBREDST}
\end{align}
where $c_{m}$ is a positive constant.
\end{lemma}
\proof When $H_{k}$ is updated by the L-BFGS formula \eqref{ILBFGS}, from Lemma
\ref{LEMHLB}, we know that $H_{k}$ is symmetric positive definite and its
eigenvalues are greater than 1/2. When $H_{k} = \left(\nabla^{2} f(x_{k})\right)^{-1}$, 
from the assumption \eqref{ASCONV} of $\nabla^{2} f(x_{k})$, we know that 
$H_{k}$  is symmetric positive definite and its eigenvalues are greater than
$1/M$. By combining these two cases, we know that the eigenvalues of $H_{k}$ 
are greater than $c_{m} = \min\{1/2, \, 1/M\}$.

\vskip 2mm

By using the eigenvalue decomposition of $H_{k}$, from the explicit continuation
method \eqref{EXPTC}-\eqref{XK1} and the quadratic model \eqref{LOAM}, we have
\begin{align}
      m_{k}(0) - m_{k}(s_{k}) \ge -\frac{1}{2}g_{k}^{T}s_{k}
       = \frac{\Delta t_{k}}{2(1+\Delta t_{k})}g_{k}^{T}H_{k}g_{k}
      = \frac{c_{m}\Delta t_{k}}{2(1+\Delta t_{k})}\|g_{k}\|^{2}.  \label{LMOLB}
\end{align}
In the first inequality in equation \eqref{LMOLB}, we use the property
$(1+0.5 \Delta t_{k})/(1+\Delta t_{k}) \ge 0.5$ when $\Delta t_{k} \ge 0$.
Consequently, we prove the result \eqref{PLBREDST}. \qed

\vskip 2mm

In order to prove that $p_{g_k}$ converges to zero when $k$ tends to infinity,
we need to estimate the lower bound of time-stepping sizes 
$\Delta t_{k} \, (k = 1, \, 2, \, \ldots)$.

\vskip 2mm

\begin{lemma} \label{DTBOUND}
Assume that $f: \; \Re^{n} \rightarrow \Re$ is twice continuously differentiable
and its gradient $g(x)$ satisfies the following Lipschitz continuity:
\begin{align}
     \|g(x) - g(y)\| \le L_{C} \|x - y\|, \; \forall x, \, y  \in S_{f}, 
     \label{LIPSCHCON}
\end{align}
where $L_{C}$ is the Lipschitz constant. The Hessian matrix $\nabla^{2} f(x)$
satisfies the strong convexity \eqref{ASCONV}. We suppose that the sequence $\{x_{k}\}$ is
generated by Algorithm \ref{ALGEPTC}. Then, there exists a positive constant
$\delta_{\Delta t}$ such that
\begin{align}
    \Delta t_{k} \ge \gamma_{2} \delta_{\Delta t}
    \label{DTGEPN}
\end{align}
holds for all $k = 1, \,  2, \, \dots$, where $\Delta t_{k}$ is adaptively
adjusted by the formulas \eqref{LOAM}-\eqref{ADTK1}.
\end{lemma}

\vskip 2mm

\proof When $H_{k}$ is updated by the L-BFGS formula \eqref{ILBFGS}, from Lemma
\ref{LEMHLB}, we know that the eigenvalues of $H_{k}$ is greater than 1/2 and
it has at least $n-2$ eigenvalues of 1. When
$|s_{k-1}^{T}y_{k-1}| \ge \theta \|s_{k-1}\|^{2}$, we denote the other two
eigenvalues of $H_{k}$ as $\mu_{1}^{k}$ and $\mu_{2}^{k}$. By substituting it
into equation \eqref{ROOTQEQ}, we obtain
\begin{align}
      \mu_{1}^{k} \mu_{2}^{k} = \frac{\|y_{k-1}\|^{2}
      \|s_{k-1}\|^{2}}{(s_{k-1}^{T}y_{k-1})^{2}} \le
      \frac{\|y_{k-1}\|^{2}
      \|s_{k-1}\|^{2}}{\theta^{2} \|s_{k-1}\|^{4}}
      = \frac{1}{\theta^{2}}\frac{\|y_{k-1}\|^{2}}{\|s_{k-1}\|^{2}}.
      \label{TWOEIGS}
\end{align}

\vskip 2mm

From Lemma \ref{LBSOAM} and Algorithm \ref{ALGEPTC}, we know $f(x_{k}) \le
f(x_{0}) \, (k = 1, \, 2, \, \ldots)$. From the Lipschitz continuity
\eqref{LIPSCHCON} of $g(\cdot)$, we have
\begin{align}
     \|y_{k-1}\| = \|g(x_{k}) - g(x_{k-1})\| \le L_{C}\|x_{k} - x_{k-1}\|
     = L_{C} \|s_{k-1}\|. \label{UPBYK}
\end{align}
By substituting it into equation \eqref{TWOEIGS} and using $\mu_{i}^{k} > \frac{1}{2}
\, (i = 1, \, 2)$, we obtain
\begin{align}
     \frac{1}{2} < \mu_{i}^{k} < \frac{2L_{C}^{2}}{\theta^{2}}, \; i = 1, \, 2.
     \label{BOUNEIG}
\end{align}
That is to say, the eigenvalues of $H_{k}$ are less than or equal to
$\max\{1, \, 2L_{C}^{2}/\theta^{2}\}$.

\vskip 2mm

When $H_{k} = \left(\nabla^{2} f(x_{k})\right)^{-1}$, from the assumption 
\eqref{ASCONV} of $\nabla^{2} f(x_{k})$, we know that $H_{k}$ is symmetric
 positive definite and its eigenvalues are less than or equal to $1/m$.
 
\vskip 2mm 

By combining these two cases, we obtain that the eigenvalues of $H_{k}$ are 
less than or equal to $M_{H}$ when $H_{k}$ is updated by the preconditioning 
formula \eqref{SHBFGS}, where 
$M_{H} = \max \{1, \, {2L_{C}^{2}}/{\theta^{2}}, \, 1/m\}$. By using the
property of the matrix norm, we have
\begin{align}
      \|H_{k}g_{k}\| \le M_{H}\|g_{k}\|.  \label{UBHPGK}
\end{align}

\vskip 2mm

From the first-order Taylor expansion, we have
\begin{align}
    f(x_{k}+ s_{k}) =  f(x_{k}) +
    \int_{0}^{1} s_{k}^{T}g (x_{k} + t s_{k}) dt.\label{FOTEFK}
\end{align}
Thus, from equations \eqref{LOAM}-\eqref{MRHOK}, \eqref{PLBREDST}, \eqref{FOTEFK}
and the Lipschitz continuity \eqref{LIPSCHCON} of $g(\cdot)$, we have
\begin{align}
      & \left|\rho_{k} - 1\right| =  \left|\frac{(f(x_{k}) - f(x_{k}+s_{k}))
       - (m_{k}(0) - m_{k}(s_{k}))}{m_{k}(0) - m_{k}(s_{k})}\right|
       \nonumber \\
      & \quad \le 
      \frac{\left|\int_{0}^{1}s_{k}^{T}(g(x_{k} + t s_{k}) - g(x_{k}))dt\right|}
      {m_{k}(0) - m_{k}(s_{k})} + \frac{0.5\Delta t_{k}}{1+0.5\Delta t_{k}} \nonumber \\
       & \quad \le  \frac{L_{C}(1+ \Delta t_{k})}{c_{m}\Delta t_{k}}
       \frac{\|s_{k}\|^{2}}{\|g_{k}\|^{2}} + \frac{0.5\Delta t_{k}}{1+0.5\Delta t_{k}}.
       \label{ESTRHOK}
\end{align}
By substituting equation \eqref{EXPTC} and equation \eqref{UBHPGK} into equation
\eqref{ESTRHOK}, we have
\begin{align}
      & \left|\rho_{k} - 1\right|
      \le \frac{L_{C}\Delta t_{k}}{c_{m}(1+\Delta t_{k})}
       \frac{\|H_{k}g_{k}\|^{2}}{\|g_{k}\|^{2}}
       + \frac{0.5\Delta t_{k}}{1+0.5\Delta t_{k}}       \nonumber \\
      & \le \frac{L_{C}\Delta t_{k}}{c_{m}(1+\Delta t_{k})}
       \frac{M_{H}^{2}\|g_{k}\|^{2}}{\|g_{k}\|^{2}}
       + \frac{0.5\Delta t_{k}}{1+0.5\Delta t_{k}}
      \le \frac{(L_{C} M_{H}^{2}+0.5c_{m})\Delta t_{k}}{c_{m}(1+0.5\Delta t_{k})}.
      \label{UBROHK}
\end{align}
We denote
\begin{align}
     \delta_{\Delta t} \triangleq \frac{c_{m}\eta_{1}}{L_{C} M_{H}^{2}+0.5c_{m}}.
     \label{UPBMPD}
\end{align}
Then, from equation \eqref{UBROHK}-\eqref{UPBMPD}, when
$\Delta t_{k} \le \delta_{\Delta t}$, it is not difficult to verify
\begin{align}
    \left|\rho_{k} - 1\right| \le \frac{L_{C} M_{H}^{2}+0.5c_{m}}{c_{m}}\Delta t_{k} \le \eta_{1}.
    \label{RHOLETA1}
\end{align}

\vskip 2mm

We assume that $K$ is the first index such that $\Delta t_{K} \le
\delta_{\Delta t}$ where $\delta_{\Delta t}$ is defined by equation \eqref{UPBMPD}.
Then, from equations \eqref{UPBMPD}-\eqref{RHOLETA1}, we know that
$|\rho_{K} - 1 | \le \eta_{1}$. According to the time-stepping adjustment
formula \eqref{ADTK1}, $x_{K} + s_{K}$ will be accepted and the time-stepping size
$\Delta t_{K+1}$ will be enlarged. Consequently, the time-stepping size $\Delta t_{k}$ holds
$\Delta t_{k}\ge \gamma_{2}\delta_{\Delta t}$ for all $k = 1, \, 2, \ldots$. \qed

\vskip 2mm

By using the results of Lemma \ref{LBSOAM} and Lemma \ref{DTBOUND}, we prove
the global convergence of Algorithm \ref{ALGEPTC} for the  unconstrained
optimization problem \eqref{UNOPT} as follows.

\vskip 2mm

\begin{theorem}
Assume that $f: \; \Re^{n} \rightarrow \Re$ is twice continuously differentiable
and its gradient $\nabla f(\cdot)$ satisfies the Lipschitz continuity \eqref{LIPSCHCON}.
The Hessian matrix $\nabla^{2} f(x)$ satisfies the strong convexity \eqref{ASCONV}.
Furthermore, we suppose that $f(x)$ is lower bounded when $x \in S_{f}$, where
the level set $S_{f}$ is defined by equation \eqref{LSCONFBD}. The sequence
$\{x_{k}\}$ is generated by Algorithm \ref{ALGEPTC}. Then, we have
\begin{align}
   \lim_{k \to \infty} \inf \|g_{k}\| = 0. \label{LIMPGKZ}
\end{align}
\end{theorem}
\proof According to Lemma  \ref{DTBOUND} and Algorithm \ref{ALGEPTC}, we know
that there exists an infinite subsequence $\{x_{k_{i}}\}$ such that the trial 
steps $s_{k_i}$ are accepted, i.e., $\rho_{k_{i}} \ge \eta_{a}, \, i=1,\, 2,\ldots$. Otherwise,
all steps are rejected after a given iteration index, then the time-stepping
size will keep decreasing, which contradicts \eqref{DTGEPN}. Therefore,
from equations \eqref{MRHOK} and \eqref{PLBREDST}, we have
\begin{align}
     & f(x_{0}) - \lim_{k \to \infty} f(x_{k})
     = \sum_{k = 0}^{\infty} (f(x_{k}) - f(x_{k+1}))
      \nonumber \\
     & \ge \eta_{a} \sum_{i = 0}^{\infty}
    \left(m_{k_{i}}(0) - m_{k_{i}}(s_{k_{i}})\right)
    \ge \eta_{a} \sum_{i = 0}^{\infty}
    \frac{c_{m}\Delta t_{k_{i}}}{2(\Delta t_{k_{i}}+1)} \|g_{k_{i}}\|.
    \label{LIMSUMFK}
\end{align}

\vskip 2mm

From the result \eqref{DTGEPN} of Lemma \ref{DTBOUND}, we know that
$\Delta t_{k} \ge \gamma_{2} \delta_{\Delta t} \, ( k = 1, \, 2, \, \dots)$.
By substituting it into equation \eqref{LIMSUMFK}, we have
\begin{align}
      f(x_{0}) - \lim_{k \to \infty} f(x_{k}) \ge \eta_{a} \sum_{i = 0}^{\infty}
    \frac{\gamma_{2} c_{m}\delta_{\Delta t}}
    {2(\gamma_{2} \delta_{\Delta t}+1)}\|g_{k_{i}}\|. \label{OBJDMD}
\end{align}
Since $f(x)$ is lower bounded when $x \in S_{f}$ and the sequence $\{f(x_{k})\}$
is monotonically decreasing, we have $\lim_{k \to \infty} f(x_{k}) = f^{\ast}$.
By substituting it into equation \eqref{OBJDMD}, we obtain the result
\eqref{LIMPGKZ}. \qed

\section{Numerical experiments}

\vskip 2mm

In this section, some numerical experiments are performed to test the performance
of Algorithm \ref{ALGEPTC} (the Eptctr method). The codes are executed by a HP
notebook with the Intel quad-core CPU and 8GB RAM. We compare Eptctr with
the trust-region method and the line searh method (the built-in subroutine
fminunc.m of the MATLAB2019a environment)
\cite{Broyden1970,CL1994,CL1996,Davidon1991,FP1963,MATLAB}
for $47$ unconstrained optimization problems which can be found in
\cite{AD2005,ML2004,SB2013}. The trust-region method and the line search method
are two classical methods for solving unconstrained optimization problems and
these two methods have been widely used until today. Therefore, we select these
two typical methods as the basis for comparison. The termination conditions of
the three compared methods are all set by
\begin{align}
    & \|\nabla f(x_{k})\|_{\infty} \le 1.0 \times 10^{-6}.
    \label{GRADTOL}
\end{align}
The initial points are set to $x_{0} = 2 \times \text{ones}(n,\, 1)$ for all test problems.

\vskip 2mm

\begin{figure}[!htbp]
\centering
\includegraphics[width=0.80\textwidth]{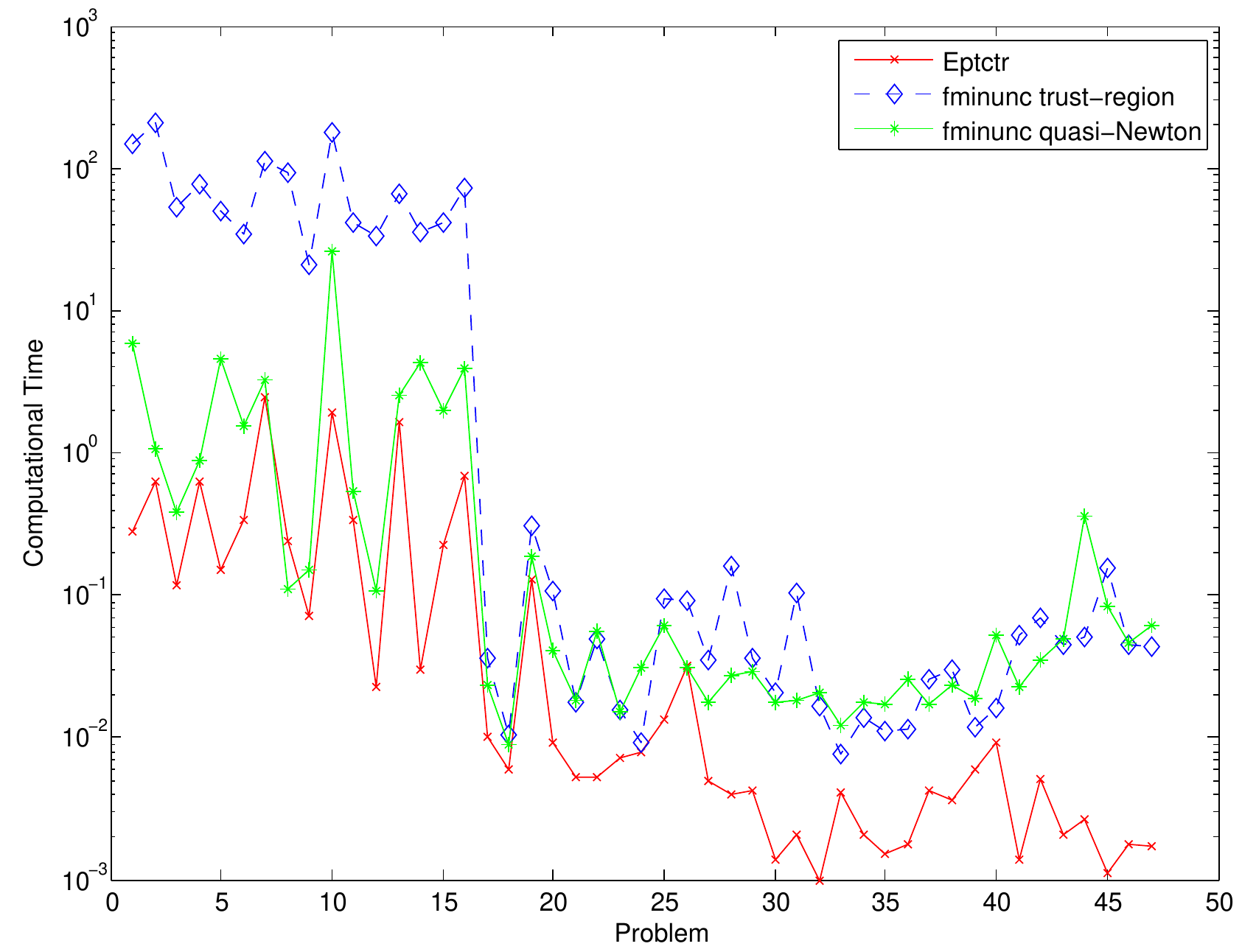}
\caption{The computational time of each problem.}\label{FIGCPU}
\end{figure}

The numerical results are arranged in Table \ref{TABCOMT1}, Table \ref{TABCOMT2}
and Figure \ref{FIGCPU}. From Tables \ref{TABCOMT1}-\ref{TABCOMT2}, we find
that Eptctr works well for those 47 test problems. However, the trust-region method
(fminunc.m is set by the trust-region method) and the line search method
(fminunc.m is set by the quasi-Newton method) fail to solve 6 problems and 13
problems, respectively.  Thus, Eptctr is more robust than the traditional optimization
method such as the trust-region method and the line search method.

\vskip 2mm

Moreover, from Table \ref{TABCOMT1}, we find that the computational time of Eptctr 
is about one percent of that of the trust-region method (fminunc.m is set by the 
trust-region method) or one fifth of that the line search method (fminunc.m is 
set by the quasi-Newton method). One of the reasons is that the generalized 
gradient flow is non-stiff in the transient-state phase and Eptctr uses 
the L-BFGS method \ref{ILBFGS} as the preconditioning technique to follow 
their trajectories well for the most problems. Consequently, Eptctr only 
involves three pairs of the inner product of two vectors to obtain the trial 
step $s_{k}$ at every iteration of the transient-state phase. 

\vskip 2mm 

However, the trust-region method needs at least to solve a linear system of 
equations and involves about $1/3 \, n^{3}$ flops (p. 169, \cite{GV2013}) 
at every iteration. Actually, in order to obtain the search step $s_{k}$, 
the trust-region method needs to solve a trust-region subproblem, which 
requires to solve a nonlinear system of equations at every iteration. The line 
search method involves one matrix-vector product to obtain the search direction 
$d_{k}$ and it involves about $n^{2}$ flops at every iteration. The other reason 
is that the step size $\alpha_{k}$ of the line search method is tried from 1 and 
reduced by the half with many times at every iteration. Since the linear model 
$f(x_{k}) + g_{k}^{T}s_{k}$ may not approximate $f(x_{k}+s_{k})$
well in the transient-state phase, the step size $\alpha_{k}$ will be small.
Consequently, the line search strategy consumes the unnecessary trial steps in
the transient-state phase.

\vskip 2mm

\begin{table}[!http]
  \newcommand{\tabincell}[2]{\begin{tabular}{@{}#1@{}}#2\end{tabular}}
  \scriptsize
  \centering
  \caption{Numerical results of large-scale problems with $n = 1000$.}
  \label{TABCOMT1}
  \resizebox{\textwidth}{!}{
  \begin{tabular}{|c|c|c|c|c|c|c|c|c|c|}
  \hline
  \multirow{2}{*}{Problems } & \multicolumn{2}{c|}{Eptctr} & \multicolumn{2}{c|}{fminunc.m (trust-region)} & \multicolumn{2}{c|}{fminunc.m (quasi-Newton)}  \\ \cline{2-7}
  & \tabincell{c}{Iter (time (s))}     & $||g(x^{it})||_\infty$  & \tabincell{c}{Iter (time (s))} & $||g(x^{it})||_\infty$   & \tabincell{c}{Iter (time (s))}  & $||g(x^{it})||_\infty$  \\ \hline

  \tabincell{c}{1. Trid Function \\ (n = 1000)}   & \tabincell{c}{28 (0.2800)} & 7.7882E-07 & \tabincell{c}{21 (147.3268)
  \\ \red{(failed)}} & \red{0.0052} &  \tabincell{c} {89 (5.8060) \\ \red{(failed)}}   &  \red{185.3336}  \\ \hline

  \tabincell{c}{2. Rosenbrock Function \\ (n = 1000)}   & \tabincell{c}{37 (0.6280)} &1.8925E-07 & \tabincell{c}{15 (208.36)} &6.9663E-11 & \tabincell{c}{23 (1.0736)} & \red{3.7971E-05}   \\ \hline

  \tabincell{c}{3. Ackley Function \\ (n = 1000)}   & \tabincell{c}{74 (0.1168)} & 2.1729E-06 & \tabincell{c}{2 (52.7352)} &6.9663E-11 & \tabincell{c}{4 (0.3761)} & 0 \\ \hline

  \tabincell{c}{4. Dixon Price Function \\ (n = 1000)}   & \tabincell{c}{37 (0.6220)} &1.8582E-07 & \tabincell{c}{11 (78.2305)
  \\ \red{(failed)}} & \red{5.1937E-04} & \tabincell{c}{95 (4.5570)\\ \red{(failed)}} & \red{10.3910}   \\ \hline

  \tabincell{c}{5. Levy Function \\ (n = 1000)}   & \tabincell{c}{38 (0.1521)} & 7.6318E-07 & \tabincell{c}{6 (49.5847)} & 4.9983E-13  & \tabincell{c}{11 (0.8782)} & 2.5648E-07 \\ \hline

  \tabincell{c}{6. Molecular Energy \\ Function (n = 1000)}   & \tabincell{c}{22 (0.3407)} &6.5946E-07  & \tabincell{c}{4 (34.4602)} & 3.1098E-09 & \tabincell{c}{6 (1.5520)} &3.6708E-06\\ \hline

  \tabincell{c}{7. Powell Function \\ (n = 1000)}   & \tabincell{c}{44 (2.4865)} &6.2502E-07 &\tabincell{c}{16 (110.4723)} &3.8521E-07 & \tabincell{c}{25 (3.2259) \\ \red{(failed)}} & \red{1.7934E-04} \\ \hline

  \tabincell{c}{8. Quartic With Noise \\ Function (n = 1000)}   & \tabincell{c}{383 (0.2384)} &3.1145E-06 & \tabincell{c}{13 (91.8643)} &5.4263E-07 & \tabincell{c}{1 (0.1120)} &0 \\ \hline

  \tabincell{c}{9. Rastrigin Function \\ (n = 1000)}   & \tabincell{c}{9 (0.0712)} &3.5796E-08 & \tabincell{c}{1 (21.1994)} & 3.5728E-09 & \tabincell{c}{2 (0.1493)} & 0\\ \hline

  \tabincell{c}{10.Rotated Hyper Ellipsoid \\ Function (n = 1000)}   & \tabincell{c}{28 (1.9075)} & 4.3487E-08 & \tabincell{c}{2 (176.9553)} &4.3032E-13 & \tabincell{c}{99 (25.9046)\\ \red{(failed)}} & \red{0.2211} \\ \hline

  \tabincell{c}{11.Schwefel Function \\ (n = 1000)}   & \tabincell{c}{36 (0.3373)} & 9.4745E-07 & \tabincell{c}{5 (40.8946)} &9.0242E-07 & \tabincell{c}{5 (0.5322)} & 0 \\ \hline

  \tabincell{c}{12.Sphere Function \\ (n = 1000)}   & \tabincell{c}{14 (0.0230)} & 8.4047E-07 & \tabincell{c}{4 (33.5744)} & 2.9759E-08 & \tabincell{c}{2 (0.1065)} & 1.4901E-08 \\ \hline

  \tabincell{c}{13. Styblinski Tang \\ Function (n = 1000)}   & \tabincell{c}{68 (1.6197)} & 2.6205E-06 & \tabincell{c}{8 (65.4069)} & 1.0887E-06 & \tabincell{c}{7 (2.5245) \\ \red{(failed)}} & \red{1.7776E-04} \\ \hline

  \tabincell{c}{14. Sum Squares \\ Function (n = 1000)}   & \tabincell{c}{22 (0.0300)} & 7.9787E-07 & \tabincell{c}{4 (35.9020)} & 1.0839E-20 & \tabincell{c}{99 (4.3612) \\ \red{(failed)}} & \red{0.4370} \\ \hline

  \tabincell{c}{15. Schubert Function \\ (n = 1000)}   & \tabincell{c}{16 (0.2256)} & 8.2217E-07 & \tabincell{c}{5 (41.7283)} & 2.1316E-14 & \tabincell{c}{6 (1.9671)} & \red{5.8561E-05} \\ \hline

  \tabincell{c}{16. Stretched V Function \\ (n = 1000)}   & \tabincell{c}{14 (0.6951)} & 4.2645E-07 & \tabincell{c}{3 (73.4679)} & 4.8573E-09 & \tabincell{c}{4 (3.8702)} & 2.0675E-07 \\ \hline

\end{tabular}}
\end{table}

\vskip 2mm

\begin{table}[!http]
  \newcommand{\tabincell}[2]{\begin{tabular}{@{}#1@{}}#2\end{tabular}}
  \scriptsize
  \centering
  \caption{Numerical results of small-scale problems with $n \leq 10$.}
  \label{TABCOMT2}
  \resizebox{\textwidth}{!}{
  \begin{tabular}{|c|c|c|c|c|c|c|c|c|c|}
  \hline
  \multirow{2}{*}{Problems } & \multicolumn{2}{c|}{Eptctr} & \multicolumn{2}{c|}{fminunc.m (trust-region)} & \multicolumn{2}{c|}{fminunc.m (quasi-Newton)}  \\ \cline{2-7}
  & \tabincell{c}{Iter (time (s))}     & $||g(x^{it})||_\infty$  & \tabincell{c}{Iter (time (s))}     & $||g(x^{it})||_\infty$   & \tabincell{c}{Iter (time)}  & $||g(x^{it})||_\infty$  \\ \hline

  \tabincell{c}{17. Beale Function\\ (n = 2)}   & \tabincell{c}{94 (0.0100)} & 3.4136E-08 & \tabincell{c}{17 (0.0359)} & \red{2.3924E-05} & \tabincell{c}{21 (0.0232)} & \red{7.1145E-05}  \\ \hline

  \tabincell{c}{18. Booth Function \\ (n = 2)}   & \tabincell{c}{24 (0.0060)} &4.1473E-07 & \tabincell{c}{1 (0.0105)} & 0 & \tabincell{c}{1 (0.0090)} & 1.2666E-06   \\ \hline

  \tabincell{c}{19. Branin Function \\ (n = 2)}   & \tabincell{c}{32 (0.1275)} & 1.8526E-07 & \tabincell{c}{4 (0.3062)} & 9.9309E-08 & \tabincell{c}{7 (0.1863)} & 4.0602E-06 \\ \hline

  \tabincell{c}{20. Easom Function \\ (n = 2)}   & \tabincell{c}{27 (0.0092)} & 9.9323E-07 & \tabincell{c}{6 (0.1055)} & 2.6839E-09 & \tabincell{c}{3 (0.0405)} & 4.7432E-09   \\ \hline

  \tabincell{c}{21. Griewank Function \\ (n = 10)}   & \tabincell{c}{28 (0.0053)} & 8.1615E-07 & \tabincell{c}{5 (0.0176)} & 7.2012E-10  & \tabincell{c}{13 (0.0180)} & 2.9802E-07 \\ \hline

  \tabincell{c}{22. Hosaki Function \\ (n = 2)}   & \tabincell{c}{25 (0.0053)} & 8.2774E-07  & \tabincell{c}{1 (0.0489)} & 8.5432E-09 & \tabincell{c}{5 (0.0552)} &3.7253E-07\\ \hline

  \tabincell{c}{23. Levy13 Function \\ (n = 2)}   & \tabincell{c}{22 (0.0071)} & 8.9312E-07 &\tabincell{c}{3 (0.0158)} & 1.4647E-11 & \tabincell{c}{2 (0.0153)} & 7.8931E-07 \\ \hline

  \tabincell{c}{24. Matyas Function \\ (n = 2)}   & \tabincell{c}{48 (0.0080)} & 7.5460E-07 & \tabincell{c}{1 (0.0092)} & 4.7684E-08 & \tabincell{c}{2 (0.0310)} & 2.1160E-07 \\ \hline

  \tabincell{c}{25. Mccormick Function \\ (n = 2)}   & \tabincell{c}{24 (0.0132)} & 4.2654E-07 & \tabincell{c}{4 (0.0949)
  \\ \red{(failed)}} & \red{6.0400E-04} & \tabincell{c}{6 (0.0617)} & 1.8692E-08\\ \hline

  \tabincell{c}{26. Perm Function \\ (n = 4)}   & \tabincell{c}{174 (0.0320)} & 2.5599E-08 & \tabincell{c}{20 (0.0908) \\ \red{(failed)}} & \red{0.0112} & \tabincell{c}{34 (0.0307)\\ \red{(failed)}} & \red{0.0503} \\ \hline

  \tabincell{c}{27. Power Sum Function \\ (n = 4)}   & \tabincell{c}{28 (0.0050)} & 2.0605E-07 & \tabincell{c}{6 (0.0353)} & 2.2251E-07 & \tabincell{c}{6 (0.0176) \\ \red{(failed)}} & \red{1.2351E-04} \\ \hline

  \tabincell{c}{28. Price Function \\ (n = 2)}   & \tabincell{c}{61 (0.0040)} & 4.4811E-07 & \tabincell{c}{11 (0.1592)} & 7.8428E-06 & \tabincell{c}{27 (0.0271) \\ \red{failed)}} & \red{9.0677E-04} \\ \hline

  \tabincell{c}{29. Zakharov Function \\ (n = 10)}   & \tabincell{c}{48 (0.0042)} & 2.7717E-07 & \tabincell{c}{15 (0.0363)} & 5.9171E-06 & \tabincell{c}{26 (0.0291)\\ \red{(failed)}} & \red{3.3111} \\ \hline

  \tabincell{c}{30. Bohachevsky Function \\ (n = 2)}   & \tabincell{c}{32 (0.0014)} & 7.8927E-07 & \tabincell{c}{9 (0.0204)} & 2.1823E-06 & \tabincell{c}{8 (0.0174)} & 9.9838E-07 \\ \hline

  \tabincell{c}{31. Colville Function \\ (n = 4)}   & \tabincell{c}{37 (0.0021)} & 8.0037E-08 & \tabincell{c}{70 (0.1032)\\ \red{(failed)}} & \red{0.0047} & \tabincell{c}{27 (0.0185)\\ \red{(failed)}} & \red{0.0014} \\ \hline

  \tabincell{c}{32. Drop Wave \\ Function (n = 2)}   & \tabincell{c}{12 (0.0010)} & 3.6227E-07 & \tabincell{c}{7 (0.0168)} & 2.7057E-08 & \tabincell{c}{4 (0.0206) } & 1.2368E-06 \\ \hline

  \tabincell{c}{33. Schaffer Function\\ (n = 2)}   & \tabincell{c}{13 (0.0041)} & 3.3965E-07 & \tabincell{c}{3 (0.0076)} & 4.0569E-10 & \tabincell{c}{4 (0.0120)} & 3.3573E-09 \\ \hline

  \tabincell{c}{34. Six Hump Camel\\Function (n = 2)}   & \tabincell{c}{32 (0.0021)} & 1.7971E-07 & \tabincell{c}{7 (0.0139)} & 2.6754E-06 & \tabincell{c}{13 (0.0176)} & \red{7.5437E-05} \\ \hline

  \tabincell{c}{35. Three Hump Camel\\Function (n = 2)}   & \tabincell{c}{28 (0.0015)} & 2.9007E-07 & \tabincell{c}{5 (0.0110)} & 1.0843E-10 & \tabincell{c}{10 (0.0172)} & 5.8845E-06 \\ \hline

  \tabincell{c}{36. Trecanni Function\\ (n = 2)}   & \tabincell{c}{27 (0.0018)} & 5.9765E-07 & \tabincell{c}{7 (0.0116)} & 1.2378E-10 & \tabincell{c}{18 (0.0259)} & 2.1014E-06 \\ \hline

  \tabincell{c}{37. Box Betts Exponential\\Quadratic Sum Function (n = 3)}   & \tabincell{c}{12  (0.0042)} & 8.5907E-08 & \tabincell{c}{7 (0.0255)} & 8.7407E-12 & \tabincell{c}{14 (0.0173)} & 1.4076E-07 \\ \hline

  \tabincell{c}{38. Chichinadz Function\\(n = 2)}   & \tabincell{c}{33  (0.0036)} & 2.1310E-07 & \tabincell{c}{12 (0.0300)} & 1.4076E-13 & \tabincell{c}{11 (0.0236)} & \red{1.2013E-05} \\ \hline

  \tabincell{c}{39. Eggholder Function\\ (n = 2)}   & \tabincell{c}{78  (0.0059)} & 1.3660E-05 & \tabincell{c}{4 (0.0118)} & 8.3427E-07 & \tabincell{c}{7 (0.0188)} & 6.0934E-08 \\ \hline

  \tabincell{c}{40. Exp2 Function\\ (n = 2)}   & \tabincell{c}{172  (0.0093)} & 2.0727E-07 & \tabincell{c}{7 (0.0160)} & 3.4218E-09 & \tabincell{c}{14 (0.0530)} & 1.1368E-07 \\ \hline

  \tabincell{c}{41. Hansen Function\\ (n = 2)}   & \tabincell{c}{15  (0.0014)} & 7.1112E-07 & \tabincell{c}{7 (0.0517)} & \red{2.2678E-05} & \tabincell{c}{8 (0.0229)} & 1.1921E-06 \\ \hline

  \tabincell{c}{42. Hartmann 3-D\\dimensional Function (n = 3)}   & \tabincell{c}{54  (0.0051)} & 6.1888E-05 & \tabincell{c}{12 (0.0698)} & \red{1.6596E-05} & \tabincell{c}{15 (0.0345)} & 9.2387E-07 \\ \hline

  \tabincell{c}{43. Holder Table\\ Function (n = 2)}   & \tabincell{c}{16 (0.0021)} & 3.6270E-07 & \tabincell{c}{8 (0.0450)} & 2.6754E-09 & \tabincell{c}{4 (0.3598)\\ \red{(failed)}} & \red{7.4925E+37} \\ \hline

  \tabincell{c}{44. Michalewicz Function \\ (n = 2)}   & \tabincell{c}{31 (0.0027)} & 4.7446E-07 & \tabincell{c}{10 (0.0513)} & 1.5065E-08 & \tabincell{c}{12 (0.0498)} & 2.7057E-08 \\ \hline

  \tabincell{c}{45. Schaffer Function N.4 \\ (n = 2)}   & \tabincell{c}{2 (0.0011)} & 1.0267E-07 & 
  \tabincell{c}{46 (0.1563)} & 9.8423E-07 & \tabincell{c}{10 (0.0825)} & \red{6.0786E-05} \\ \hline

  \tabincell{c}{46. Trefethen 4 Function \\ (n = 2)}   & \tabincell{c}{29 (0.0018)} & 2.4417E-07 & \tabincell{c}{8 (0.0443)
  \\ \red{(failed)}} & \red{1.2562E-04} & \tabincell{c}{14 (0.0464) \\ \red{(failed)}} & \red{0.0052}\\ \hline

  \tabincell{c}{47. Zettl Function \\ (n = 2)}   & \tabincell{c}{40 (0.0017)} & 5.1869E-07 & \tabincell{c}{13 (0.0438)} & 3.7740E-08 & \tabincell{c}{18 (0.0608)} & 9.7767E-06 \\ \hline

\end{tabular}}
\end{table}

\vskip 2mm

\vskip 2mm

\section{Conclusions}

\vskip 2mm

For the unconstrained optimization problem, we consider an explicit continuation 
method with the trust-region updating strategy and the switching preconditioning 
technique (Eptctr) to solve it. For the well-conditioned phase,
Eptctr uses the L-BFGS method \eqref{ILBFGS} as the preconditioning technique.
Otherwise, Eptctr uses the inverse $\left(\nabla^{2} f(x_{k})\right)^{-1}$ 
of the Hessian matrix as the pre-conditioner in order to improve its robustness. 
Consequently, for the well-conditioned phase, Eptctr only involves the $6n$ flops 
to obtain its trial step $s_{k}$ at every iteration. 

\vskip 2mm 

However, the trust-region method needs to solve a trust-region subproblem, which 
requires to solve a nonlinear system of equations at every iteration. 
The line search method involves one matrix-vector product to obtain the search 
direction $d_{k}$ and it involves about $n^{2}$ flops at every iteration. 
Furthermore, the step size $\alpha_{k}$ of the line search method is tried 
from 1 and reduced by the half with many times at every iteration. Since the 
linear model $f(x_{k}) + g_{k}^{T}s_{k}$ may not approximate $f(x_{k}+s_{k})$
well in the transient-state phase, the step size $\alpha_{k}$ will be small.
Consequently, the line search strategy consumes the unnecessary trial steps in
the transient-state phase.

\vskip 2mm 

Numerical results also show that Eptctr is more robust and faster than the 
traditional optimization method such as the trust-region method and the line 
search method. The computational time of Eptctr is about one percent
of that of the trust-region method (fminunc.m is set by the trust-region method )
 or one fifth of that of the line search method (fminunc.m is set by the 
 quasi-Newton method) for the test problems.  Therefore, Eptctr can be regarded as 
 a work horse for the unconstrained 
optimization problem and it is worth exploring further. We will extend it 
to the constrained optimization problem in the future.

\section*{Acknowledgments} This work was supported in part by Grant 61876199 
from National Natural Science Foundation of China, Grant YBWL2011085 from Huawei 
Technologies Co., Ltd., and Grant YJCB2011003HI from the Innovation Research Program 
of Huawei Technologies Co., Ltd.. The authors are grateful to the anonymous 
referees for their comments and suggestions which greatly improve the presentation 
of this paper. 

\end{document}